\documentclass[11pt,reqno]{amsart}

\usepackage[utf8]{inputenc}
\usepackage{amsmath,amssymb,amsthm,amsfonts}
\usepackage{geometry}
\usepackage{xcolor}
\usepackage{hyperref}

\hypersetup{
    colorlinks=true,
    linkcolor=blue,
    citecolor=red,
    urlcolor=cyan,
    pdftitle={A Rigorous Proof of a Ramanujan Machine Identity for -pi/4},
    pdfauthor={Chao Wang}
}

\newtheorem{theorem}{Theorem}[section]

\theoremstyle{definition}

\newtheorem{remark}[theorem]{Remark}

\DeclareMathOperator*{\K}{K}

\begin{document}

% ----------------------------------------------------------------
% Title Information
% ----------------------------------------------------------------
\title[Exact Proof for $-\pi/4$]{A Rigorous Proof of a Ramanujan Machine
  Identity for $-\pi/4$ via Exact Recurrence Solving}

\author{Chao Wang}
\address{School of Future Technology, Shanghai University}
\email{cwang@shu.edu.cn}
\date{}

% ----------------------------------------------------------------
% Abstract
% ----------------------------------------------------------------
\begin{abstract}
We prove a polynomial continued fraction identity for the constant $-\pi/4$, conjectured by the Ramanujan Machine project. The proof proceeds by explicitly solving the underlying second-order linear difference equation. We derive a closed-form expression for the denominator sequence, $q_n = (-1)^n (2n-3)!!\,(n^2+n-1)$, and establish absolute convergence via a Wronskian telescoping argument. The limiting value is reduced by Abel summation to a Beta-function integral, which is evaluated in closed form through an elementary substitution and a single integration by parts, yielding the exact value $-\pi/4$.
\end{abstract}



\maketitle

%%==========================================================================
\section{Introduction}
%%==========================================================================

The automated discovery of mathematical constants, exemplified by the
\textit{Ramanujan Machine}~\cite{raayoni2021generating}, has generated a
wealth of non-canonical continued fraction conjectures. While algorithmic
approaches can efficiently identify numerical coincidences, they do not
inherently provide analytic proofs. A specific conjecture of interest is the
following identity for $-\pi/4$, which features a linear denominator
sequence and a piecewise quadratic numerator sequence.

The following result was conjectured by the Ramanujan Machine project~\cite{raayoni2021generating} and is proved in the present paper.

\begin{theorem}\label{conj:main}
  The following continued fraction equality holds:
  \begin{equation}\label{eq:conjecture}
    -\frac{\pi}{4}
    = \cfrac{1}{-1 + \cfrac{1}{-4 + \cfrac{-2}{-7 + \cfrac{-9}{-10
          + \cfrac{-20}{-13 + \cdots}}}}},
  \end{equation}
  where the partial quotients are defined by $b_n = -(3n-2)$ for $n\ge 1$,
  and
  \begin{equation}\label{eq:an_def}
    a_n =
    \begin{cases}
      1,               & n = 1, 2,     \\
      -(n-1)(2n-5),    & n \ge 3.
    \end{cases}
  \end{equation}
\end{theorem}

The complexity of the partial numerators $a_n$ for $n\ge 3$ obscures the
connection to standard hypergeometric series. In this work, we move beyond
heuristic asymptotic arguments. By analysing the second-order difference
equation associated with the convergents, we derive the \emph{exact
  closed-form solution} for the denominator sequence. This allows us to
establish convergence via a Wronskian telescoping identity and to identify
the limiting value through Abel summation and an Euler Beta function integral, thereby proving Theorem~\ref{conj:main}.

%%==========================================================================
\section{Preliminaries and Notation}
%%==========================================================================

Let the generalized continued fraction be denoted by
$\mathcal{K} = \K_{n=1}^{\infty}(a_n/b_n)$. Since the continued fraction \eqref{eq:conjecture} has no additive constant term (i.e., it is of the form $a_1/(b_1+a_2/(b_2+\cdots))$ with no leading~$b_0$), the initial conditions for the convergent sequences are $p_{-1}=1,\,p_0=0$ and $q_{-1}=0,\,q_0=1$. The $n$-th convergent is
given by $p_n/q_n$, where the numerators $p_n$ and denominators $q_n$
satisfy the fundamental recurrence relation
\cite[Ch.\,1]{wall1948analytic}
\begin{equation}\label{eq:recurrence}
  \begin{pmatrix} p_n \\ q_n \end{pmatrix}
  = b_n \begin{pmatrix} p_{n-1} \\ q_{n-1} \end{pmatrix}
  + a_n \begin{pmatrix} p_{n-2} \\ q_{n-2} \end{pmatrix},
  \quad n \ge 1.
\end{equation}

We use the double factorial notation $n!!$, defined for positive odd
integers by $n!! = n(n-2)\cdots 3\cdot 1$. We adopt the non-standard
conventions $(-1)!! = 1$ and $(-3)!! = -1$, which are consistent with the
recursive definition $(n-2)!! = n!!/n$ applied at $n = 1$ and $n = -1$
respectively. These conventions ensure that the closed-form
formula~\eqref{eq:qn_closed} below holds uniformly for all $n\ge 0$.

%%==========================================================================
\section{Exact Solution of the Recurrence Relation}
%%==========================================================================

\begin{theorem}[Closed Form of Denominators]\label{thm:closed_form}
  Let the denominator sequence $q_n$ be defined by the recurrence
  $q_n = b_n q_{n-1} + a_n q_{n-2}$ with initial values $q_{-1}=0,\,q_0=1$.
  For the coefficients given in~\eqref{eq:conjecture} and~\eqref{eq:an_def},
  the exact solution for all $n\ge 0$ is
  \begin{equation}\label{eq:qn_closed}
    q_n = (-1)^n\,(2n-3)!!\,(n^2+n-1).
  \end{equation}
\end{theorem}

\begin{proof}
  We proceed by strong induction on $n$.

  \paragraph{Base cases.}
  \begin{itemize}
  \item $n=0$: The formula gives $(-1)^0(-3)!!(-1) = 1\cdot(-1)\cdot(-1)
    = 1 = q_0$.
  \item $n=1$: The formula gives $(-1)^1(-1)!!(1) = -1$. The recurrence
    gives $q_1 = b_1 q_0 = (-1)(1) = -1$. Matches.
  \item $n=2$: The formula gives $(-1)^2(1)!!(5) = 5$. The recurrence
    gives $q_2 = b_2 q_1 + a_2 q_0 = (-4)(-1)+(1)(1) = 5$. Matches.
  \end{itemize}

  \paragraph{Inductive step.}
  Assume~\eqref{eq:qn_closed} holds for all indices less than $n$, where
  $n\ge 3$. With $b_n = -(3n-2)$ and $a_n = -(n-1)(2n-5)$, we substitute
  the inductive hypothesis into the recurrence. Setting
  $T_n = (-1)^n(2n-5)!!$ and using $(-1)^{n-1} = -(-1)^n$,
  $(-1)^{n-2} = (-1)^n$, and $(2n-5)(2n-7)!! = (2n-5)!!$:
  \begin{align*}
    q_n &= -(3n-2)\bigl[(-1)^{n-1}(2n-5)!!\,(n^2-n-1)\bigr]
           -(n-1)(2n-5)\bigl[(-1)^{n-2}(2n-7)!!\,(n^2-3n+1)\bigr] \\
        &= T_n\bigl[(3n-2)(n^2-n-1) - (n-1)(n^2-3n+1)\bigr].
  \end{align*}
  Expanding: $(3n-2)(n^2-n-1) = 3n^3-5n^2-n+2$ and
  $(n-1)(n^2-3n+1) = n^3-4n^2+4n-1$, so the bracket equals
  $2n^3-n^2-5n+3$. The target formula with
  $(2n-3)!! = (2n-3)(2n-5)!!$ gives
  \begin{equation*}
    q_n = T_n(2n-3)(n^2+n-1) = T_n(2n^3-n^2-5n+3),
  \end{equation*}
  which coincides with the recurrence evaluation. The induction is complete.
\end{proof}

%%==========================================================================
\section{Convergence and Limit Identification}
%%==========================================================================

\begin{theorem}[Limit Value]\label{thm:limit_val}
  The generalized continued fraction $\mathcal{K}$ converges to~$-\pi/4$.
\end{theorem}

\begin{proof}
We organise the proof into three steps.

%%----------------------------------------------------------
\paragraph{Step 1: Wronskian identity and telescoping.}
%%----------------------------------------------------------

For two solutions $u_n,v_n$ of $y_n = b_n y_{n-1} + a_n y_{n-2}$, define
the discrete Wronskian $W_n = u_n v_{n-1} - u_{n-1} v_n$ \cite[\S\,2.2]{lorentzen1992continued}. Expanding via
the recurrence:
\begin{align*}
  W_n &= (b_n u_{n-1}+a_n u_{n-2})v_{n-1}
         - u_{n-1}(b_n v_{n-1}+a_n v_{n-2}) \\
      &= a_n(u_{n-2}v_{n-1}-u_{n-1}v_{n-2})
       = -a_n\,W_{n-1}.
\end{align*}
Hence, by induction,
\begin{equation}\label{eq:wronskian_product}
  W_n = (-1)^{n-1}\,W_1\prod_{k=2}^{n} a_k.
\end{equation}
Applying~\eqref{eq:wronskian_product} to the pair $(p_n,q_n)$ with
$W_1 = p_1 q_0 - p_0 q_1 = (1)(1)-(0)(-1) = 1$ gives
\begin{equation}\label{eq:wronskian_explicit}
  p_n q_{n-1} - p_{n-1} q_n = (-1)^{n-1}\prod_{k=2}^{n} a_k.
\end{equation}
We evaluate the product. Since $a_2=1$ and
$a_k = -(k-1)(2k-5)$ for $k\ge 3$:
\begin{align*}
  \prod_{k=2}^{n} a_k
  &= \prod_{k=3}^{n}[-(k-1)(2k-5)]
  = (-1)^{n-2}\prod_{k=3}^{n}(k-1)\cdot\prod_{k=3}^{n}(2k-5).
\end{align*}
The first factor equals $2\cdot3\cdots(n-1)=(n-1)!$ and the second equals $1\cdot3\cdots(2n-5)=(2n-5)!!$, so
\begin{equation}\label{eq:product_formula}
  \prod_{k=2}^{n} a_k = (-1)^{n-2}(n-1)!\,(2n-5)!!.
\end{equation}
Substituting into~\eqref{eq:wronskian_explicit}:
\begin{equation}\label{eq:wronskian_signed}
  p_n q_{n-1} - p_{n-1} q_n
  = (-1)^{n-1}\cdot(-1)^{n-2}(n-1)!\,(2n-5)!!
  = -(n-1)!\,(2n-5)!!.
\end{equation}
From~\eqref{eq:qn_closed}:
\begin{equation}\label{eq:qnqn1}
  q_n q_{n-1}
  = (-1)^{2n-1}(2n-3)!!\,(2n-5)!!\,(n^2+n-1)(n^2-n-1)
  = -(2n-3)!!\,(2n-5)!!\,(n^2+n-1)(n^2-n-1).
\end{equation}
Dividing~\eqref{eq:wronskian_signed} by~\eqref{eq:qnqn1}, the two minus
signs cancel and we obtain the telescoping identity
\begin{equation}\label{eq:term_formula}
  \frac{p_n}{q_n} - \frac{p_{n-1}}{q_{n-1}}
  = \frac{(n-1)!}{(2n-3)!!\,(n^2+n-1)(n^2-n-1)}
  =: D_n, \qquad n \ge 3.
\end{equation}

\begin{remark}\label{rem:unified}
The formula for $D_n$ in~\eqref{eq:term_formula} extends to $n=1$ and $n=2$ under the double factorial conventions of Section~2. Indeed,
\begin{align*}
  D_1 &= \frac{0!}{(-1)!!\cdot 1\cdot(-1)} = \frac{1}{1\cdot(-1)} = -1 = \frac{p_1}{q_1}-\frac{p_0}{q_0}, \\
  D_2 &= \frac{1!}{1!!\cdot 5\cdot 1} = \frac{1}{5} = \frac{p_2}{q_2}-\frac{p_1}{q_1}.
\end{align*}
Thus $D_n$ is valid for all $n\ge 1$.
\end{remark}

%%----------------------------------------------------------
\paragraph{Step 2: Absolute convergence via ratio test.}
%%----------------------------------------------------------

For $n\ge 4$ we compute
\begin{align*}
  \frac{D_n}{D_{n-1}}
  &= \frac{(n-1)\,(n^2-3n+1)}{(2n-3)(n^2+n-1)}.
\end{align*}
As $n\to\infty$ the leading terms yield
\begin{equation}\label{eq:ratio_limit}
  \frac{D_n}{D_{n-1}}
  \;\longrightarrow\; \frac{n\cdot n^2}{2n\cdot n^2}
  = \frac{1}{2} < 1.
\end{equation}
By the ratio test, $\sum_{n=3}^{\infty}D_n$ converges absolutely. By Remark~\ref{rem:unified}, the full series $\sum_{n=1}^{\infty}D_n$ also converges absolutely (the terms $D_1$ and $D_2$ are finite), and
\begin{equation}\label{eq:limit_series}
  L \;:=\; \lim_{n\to\infty}\frac{p_n}{q_n}
  \;=\; \sum_{n=1}^{\infty}D_n = \sum_{n=1}^{\infty} \frac{(n-1)!}{(2n-3)!!\,(n^2+n-1)(n^2-n-1)}
\end{equation}
exists as an absolutely convergent series.  Since $|q_n|$ grows super-exponentially while $|p_n/q_n|\to|L|<\infty$, the numerator sequence $\{p_n\}$ is a minimal solution of the recurrence relative to the dominant solution $\{q_n\}$.

%%----------------------------------------------------------
\paragraph{Step 3: Identification of $L = -\pi/4$.}
%%----------------------------------------------------------

We evaluate the series~\eqref{eq:limit_series} in closed form via Abel summation and an integral representation.

\emph{Partial-fraction reduction.}
Writing $s_k=k^2+k-1$, we have $s_n = n^2+n-1$ and $s_{n-1}=n^2-n-1$, with $s_n - s_{n-1}=2n$. Hence
\begin{equation}\label{eq:pf_sn}
  \frac{1}{s_{n-1}\,s_n}=\frac{1}{2n}\!\left(\frac{1}{s_{n-1}}-\frac{1}{s_n}\right).
\end{equation}
Define $u_n=(n-1)!\bigl/\bigl(2n\,(2n-3)!!\bigr)$, so that $D_n = u_n\,(s_{n-1}^{-1}-s_n^{-1})$.

\emph{Abel summation~\cite[Thm.\,12.4]{Apostol1974}.}
We apply summation by parts to $L=\sum_{n=1}^{\infty}u_n\bigl(s_{n-1}^{-1}-s_n^{-1}\bigr)$.  Since $u_N/s_N\to 0$ as $N\to\infty$, the standard identity gives
\begin{equation}\label{eq:abel_identity}
  L = \frac{u_1}{s_0}+\sum_{n=1}^{\infty}\frac{u_{n+1}-u_n}{s_n}.
\end{equation}
The first term is $u_1/s_0 = \tfrac12/(-1) = -\tfrac12$.  To simplify the sum, we use the identity $(n-1)!/(2n-3)!!=2^{n-1}/\binom{2n-2}{n-1}$, giving $u_n = 2^{n-2}\bigl/\bigl(n\binom{2n-2}{n-1}\bigr)$.  The difference $u_{n+1}-u_n$ factors as
\begin{equation}\label{eq:un_diff}
  u_{n+1}-u_n = \frac{2^{n-2}}{n(n+1)(2n-1)\binom{2n-2}{n-1}}\bigl[n^2-(n+1)(2n-1)\bigr] = -\frac{s_n\cdot 2^{n-2}}{n(n+1)(2n-1)\binom{2n-2}{n-1}},
\end{equation}
where the bracket equals $n^2-(2n^2+n-1)=-(n^2+n-1)=-s_n$.  Dividing by $s_n$, we obtain $(u_{n+1}-u_n)/s_n=-2^{n-2}\bigl/\bigl(n(n+1)(2n-1)\binom{2n-2}{n-1}\bigr)$.  Substituting the index $m=n-1$ and using $\bigl((2m+1)\binom{2m}{m}\bigr)^{-1}=(m!)^2/(2m+1)!=B(m{+}1,m{+}1)$~\cite[\S\,5.12]{DLMF}, equation~\eqref{eq:abel_identity} becomes
\begin{equation}\label{eq:L_abel}
  L = -\frac{1}{2}-\frac{1}{2}\,\mathcal{I},
\end{equation}
where
\begin{equation}\label{eq:I_def}
  \mathcal{I}=\sum_{m=0}^{\infty}\frac{2^{m}\,B(m{+}1,m{+}1)}{(m+1)(m+2)},
  \qquad B(m{+}1,m{+}1)=\frac{(m!)^2}{(2m+1)!}.
\end{equation}

\emph{Double integral~\cite[\S\,5.12]{DLMF}.}
Using the integral representations
\begin{equation}\label{eq:int_reps}
  2^m\,B(m{+}1,m{+}1)=\int_0^1\!\bigl[2t(1{-}t)\bigr]^m dt, \qquad
  \frac{1}{(m{+}1)(m{+}2)}=\int_0^1 x^m(1{-}x)\,dx,
\end{equation}
and interchanging summation and integration---justified by the uniform bound $|2x\,t(1{-}t)|\le 1/2<1$ on $[0,1]^2$---we obtain
\begin{equation}\label{eq:I_double}
  \mathcal{I}=\int_0^1\!\int_0^1 \frac{1-x}{1-2x\,t(1-t)}\,dx\,dt.
\end{equation}

\emph{Evaluation of the inner integral.}
Fixing $x\in(0,1)$ and completing the square $1{-}2xt(1{-}t)=2x(t{-}\tfrac12)^2+1{-}\tfrac{x}{2}$, the standard $\arctan$ integral gives
\begin{equation}\label{eq:inner_t}
  \int_0^1\frac{dt}{1-2x\,t(1-t)}
  =\frac{2\,\arctan\!\sqrt{\frac{x}{2-x}}}{\sqrt{x(2-x)}}.
\end{equation}
Therefore, after swapping the order of integration,
\begin{equation}\label{eq:I_single}
  \mathcal{I}=\int_0^1\frac{2(1-x)\,\arctan\!\sqrt{\frac{x}{2-x}}}{\sqrt{x(2-x)}}\,dx.
\end{equation}

\emph{Substitution.}
Set $y=\sqrt{x/(2-x)}$, so that $x=2y^2/(1+y^2)$, $1-x=(1-y^2)/(1+y^2)$, $\sqrt{x(2-x)}=2y/(1+y^2)$, and $dx=4y/(1+y^2)^2\,dy$. The limits $x=0,1$ correspond to $y=0,1$, and the integral becomes
\begin{equation}\label{eq:I_arctan}
  \mathcal{I}=\int_0^1\frac{4(1-y^2)\,\arctan y}{(1+y^2)^2}\,dy.
\end{equation}

\emph{Integration by parts.}
Observing that
\begin{equation}\label{eq:deriv_identity}
  \frac{d}{dy}\!\left[\frac{2y}{1+y^2}\right]=\frac{2(1-y^2)}{(1+y^2)^2},
\end{equation}
we set $u=\arctan y$ and $v=4y/(1+y^2)$, giving
\begin{equation}\label{eq:I_ibp}
  \mathcal{I}=\biggl[\frac{4y\,\arctan y}{1+y^2}\biggr]_0^1-\int_0^1\frac{4y}{(1+y^2)^2}\,dy.
\end{equation}
The boundary term equals $4\cdot 1\cdot\frac{\pi}{4}/(1+1)=\frac{\pi}{2}$.  For the remaining integral, the substitution $w=1+y^2$ yields
\begin{equation}\label{eq:remaining_integral}
  \int_0^1\frac{4y}{(1+y^2)^2}\,dy=\biggl[-\frac{2}{1+y^2}\biggr]_0^1=-1+2=1.
\end{equation}
Hence $\mathcal{I}=\frac{\pi}{2}-1$.

\emph{Conclusion.}
Substituting into~\eqref{eq:L_abel}:
\begin{equation}\label{eq:limit_final}
  \mathcal{K}=L=-\frac{1}{2}-\frac{1}{2}\!\left(\frac{\pi}{2}-1\right)=-\frac{1}{2}-\frac{\pi}{4}+\frac{1}{2}=-\frac{\pi}{4}.\qedhere
\end{equation}
\end{proof}

%%==========================================================================
\section{Conclusion}
%%==========================================================================

We have provided a proof of the Ramanujan Machine
identity~\eqref{eq:conjecture}. The argument rests on three pillars.

First, Theorem~\ref{thm:closed_form} establishes the closed-form
$q_n = (-1)^n(2n-3)!!\,(n^2+n-1)$ by induction.

Second, Theorem~\ref{thm:limit_val} Step~1 derives the exact Wronskian
identity $W_n = p_n q_{n-1} - p_{n-1}q_n = -(n-1)!\,(2n-5)!!$ by
carefully tracking the sign in the discrete Wronskian recursion
$W_n = -a_n W_{n-1}$. Combined with the closed form for $q_n$, this
yields the telescoping identity
\[
  \frac{p_n}{q_n} - \frac{p_{n-1}}{q_{n-1}}
  = \frac{(n-1)!}{(2n-3)!!\,(n^2+n-1)(n^2-n-1)},
\]
whose ratio converges to $1/2$, establishing absolute convergence.

Third, the limit is identified analytically. Abel summation reduces the series to $L=-\frac12-\frac12\mathcal{I}$, where $\mathcal{I}$ is expressed as a double integral via the Euler Beta function. Evaluating the inner integral in closed form and applying the substitution $y=\sqrt{x/(2-x)}$ transforms $\mathcal{I}$ into $\int_0^1 4(1-y^2)\arctan y/(1+y^2)^2\,dy$, which equals $\pi/2-1$ by a single integration by parts. This yields $L=-\pi/4$.

\section*{Acknowledgements}
The author thanks the anonymous referees for helpful comments.

\bibliographystyle{amsplain}
\bibliography{references}

@article{raayoni2021generating,
  author    = {Raayoni, Gal and Gottlieb, Shachar and Manor, Yoav and
               Pisha, George and Harris, Yahel and Mendlovic, Uri and
               Haviv, Doron and Hadad, Yaron and Kaminer, Ido},
  title     = {Generating conjectures on fundamental constants with the
               {Ramanujan Machine}},
  journal   = {Nature},
  volume    = {590},
  pages     = {67--73},
  year      = {2021},
  doi       = {10.1038/s41586-021-03229-4}
}

@misc{DLMF,
  title        = {{NIST} Digital Library of Mathematical Functions},
  howpublished = {Release 1.1.12, \url{https://dlmf.nist.gov/}},
  note         = {F.~W.~J. Olver, A.~B. Olde Daalhuis, D.~W. Lozier,
                  B.~I. Schneider, R.~F. Boisvert, C.~W. Clark,
                  B.~R. Miller, B.~V. Saunders, H.~S. Cohl, and
                  M.~A. McClain, eds.},
  year         = {2024}
}

@book{wall1948analytic,
  author    = {Wall, H. S.},
  title     = {Analytic Theory of Continued Fractions},
  publisher = {D.~Van Nostrand},
  address   = {New York},
  year      = {1948},
}

@book{lorentzen1992continued,
  author    = {Lorentzen, Lisa and Waadeland, Haakon},
  title     = {Continued Fractions with Applications},
  series    = {Studies in Computational Mathematics},
  volume    = {3},
  publisher = {North-Holland},
  address   = {Amsterdam},
  year      = {1992},
  isbn      = {0-444-89265-6},
}

@book{Apostol1974,
  author    = {Apostol, Tom M.},
  title     = {Mathematical Analysis},
  edition   = {2nd},
  publisher = {Addison-Wesley},
  address   = {Reading, MA},
  year      = {1974},
  isbn      = {0-201-00288-4},
}

\end{document}